\documentclass[draft]{amsart}
\usepackage{amssymb}
\usepackage{latexsym}
\usepackage{amsfonts}
\usepackage{amsmath}

\newcommand{\E}{\mathcal E}

\newcommand{\cd}[2]{{\sf CD}_{\rm #1}(#2)}

\newtheorem{theorem}{Theorem}[section]

\newtheorem{lemma}[theorem]{Lemma}
\newtheorem{corollary}[theorem]{Corollary}

\theoremstyle{definition}
\newtheorem{example}[theorem]{Example}
\newcommand{\sy}[1]{{\sf S}_{#1}}

\newcommand{\soc}{\operatorname{Soc}}
\newcommand{\sym}[1]{{\sf Sym}\,#1}

\renewcommand{\wr}{\,{\sf wr}\,}
\newcommand{\hol}{{\sf Hol}\,}

\newcommand{\aut}[1]{{\sf Aut}\,{#1}}

\newcommand{\inn}[1]{{\sf Inn}\,#1}

\newcommand{\cent}[2]{{\mathbb C}_{#1}(#2)}
\newcommand{\norm}[2]{{\mathbb N}_{#1}\left(#2\right)}

\renewcommand{\leq}{\leqslant}
\renewcommand{\geq}{\geqslant}

\font\tenvr=cmmi10 scaled 1600
\textfont"F=\tenvr
\renewcommand{\wr}
    {\mathrel{\mkern-1mu\mathchar"0F7B\mkern-1mu}}

\title[Quasiprimitive groups and blow-up decompositions]{Quasiprimitive 
groups and blow-up decompositions}
\author{Robert W. Baddeley, Cheryl E. Praeger, Csaba Schneider}
\address[Baddeley]{
32 Arbury Road\\Cambridge CB4 2JE, UK\\
\tt robert.baddeley@ntworld.com}

\address[Praeger]{School of Mathematics and Statistics\\
The University of Western Australia\\
35 Stirling Highway 6009 Crawley\\
Western Australia\\
\tt praeger@maths.uwa.edu.au\\
www.maths.uwa.edu.au/$\sim$praeger}
\address[Schneider]{Informatics Laboratory\\ 
Computer and Automation Research Institute\\ 
The Hungarian Academy of Sciences\\ 
1518 Budapest Pf.\ 63\\
Hungary\\
\tt csaba.schneider@sztaki.hu\\
www.sztaki.hu/$\sim$schneider}

\begin{document}

\date{draft typeset \today}
\thanks{The authors acknowledge the support of an Australian Research Council 
Discovery grant.}

\begin{abstract}
The blow-up construction by L.\ G.\ Kov\'acs has been a very useful tool 
to study embeddings of finite primitive permutation groups into wreath products in product 
action. In the present paper we extend the concept of a blow-up to 
finite quasiprimitive permutation
groups, and use it to study embeddings of finite
quasiprimitive groups into wreath products.
\end{abstract}
\maketitle

\section{Introduction}

It is an important problem in the study of permutation groups to describe,
in as much detail as possible, the  inclusions among different 
classes of groups. For primitive groups this problem was solved 
by the second author in~\cite{prae:inc}. The special case of describing
the possible inclusions of primitive groups into wreath products in product 
action relied on the concept of a blow-up defined by L.\ G.\ Kov\'acs
in his seminal paper~\cite{blowups}. 

A finite permutation group
is said to be {\em quasiprimitive} if
all its non-trivial normal subgroups are transitive.
In our research into quasiprimitive permutation groups and their
actions on combinatorial objects we found it necessary to extend the 
results of~\cite{prae:inc} and to study the class of inclusions of
quasiprimitive groups into wreath products in product action.
It turned out that the class of these inclusions is much richer 
than that of the primitive groups. In order to give a detailed description
in our case, we generalised in \cite{recog} the concept of a system of 
product imprimitivity introduced by Kov\'acs in~\cite{decomp} and 
defined the concept of a Cartesian decomposition.
The stabiliser, in a finite symmetric group, of a homogeneous 
Cartesian decomposition
of the underlying set is a wreath product, in product action, of 
smaller symmetric groups. Thus the problem of finding the set of such wreath products in product action 
that contain a given permutation group $G$ is equivalent to finding
all $G$-invariant homogeneous Cartesian 
decompositions of the underlying set. The details of this work can be found 
in~\cite{recog,transcs,intrans,3types}. Our efforts made it possible 
in~\cite{bad:quasi} 
to describe satisfactorily the inclusions of quasiprimitive groups into wreath 
products in product action. The philosophy behind this description is the same
as  in~\cite{prae:inc}: define several classes of natural inclusions, and
give a complete description of the exceptional inclusions. 

The natural inclusions in the case of 
primitive groups correspond to blow-up decompositions as defined by 
Kov\'acs~\cite{blowups}. We found it necessary to define a similar concept
for quasiprimitive groups. However, we also want this new concept to fit 
into our more general combinatorial framework. 
Hence, we first introduce the class
of normal Cartesian decompositions. Informally speaking, if $\E=\{\Gamma_1,
\ldots,\Gamma_\ell\}$ 
is a normal Cartesian
decomposition for a permutation group $G$, then $G$ has a transitive 
normal subgroup $M$
which can be written as a direct product $M=M_1\times\cdots\times M_\ell$
in such a way that $M_1\times\cdots\times M_\ell$ acts in product action on $\Gamma_1\times\cdots\times\Gamma_\ell$. 
We define a blow-up
decomposition for a permutation group as a special type of normal 
Cartesian decomposition.
It is not hard to see that 
the blow-up decomposition given by Kov\'acs in~\cite{blowups} is equivalent 
to our blow-up decomposition defined in Section~\ref{sec2} for the
class of primitive groups whose socle is non-regular; see also Theorem~\ref{r2.5}.

The aim of this paper is to study blow-up decompositions of 
quasiprimitive groups. In particular we investigate the 
extent to which a  quasiprimitive
permutation group can be recovered from knowledge of 
its components under a blow-up 
decomposition. We also prove an important theorem that was stated without proof, and used, 
in~\cite{bad:quasi} (see Theorem~\ref{main}). 

The structure of the paper is as follows. In Section~\ref{newsec}
we review some well-known facts about wreath products and quasiprimitive groups. In Section~\ref{sec2} we 
give the definition of a Cartesian decomposition, and we define normal 
decompositions and blow-up decompositions. In Section~\ref{sec3} 
we investigate the relationship between a quasiprimitive group and its components
under a blow-up decomposition. In Section~\ref{sec4} we study normal
decompositions that are not blow-ups. Finally in Section~\ref{sec5} we
state and 
prove Theorem~\ref{main} which was used in the study~\cite{bad:quasi} of inclusions
of quasiprimitive groups. 

In this 
paper we use the following notation. Permutations act on the right: if $\pi$ is 
a permutation and
$\omega$ is a point then the image of $\omega$ under $\pi$ is denoted
$\omega\pi$.  If $G$ is a group acting on a set $\Omega$ and $\Gamma$ is a
subset of $\Omega$, then $G_\Gamma$ and $G_{(\Gamma)}$ denote respectively
the setwise and the pointwise  stabiliser in $G$ of $\Gamma$.
All groups that appear in this paper are finite.

\section{Primitive and quasiprimitive groups}
\label{newsec}

First in this section we review wreath products and their product actions,
as they play an important part in our research.
Let $\Gamma$ be a finite set, $L\leq\sym\Gamma$, $\ell\geq 2$ an
integer, and $H\leq\sy\ell$. The {\em
wreath product} $L\wr H$ is the semidirect product
$L^\ell\rtimes H$, where, for
$(x_1,\ldots,x_\ell)\in L^\ell$ and $\sigma\in H$,
$(x_1,\ldots,x_\ell)^{\sigma^{-1}}=(x_{1{\sigma}},\ldots,x_{\ell{\sigma}})$. The {\em product
action} of $L\wr H$ is the action of $L\wr H$ on $\Gamma^\ell$ defined
by
$$
(\gamma_1,\ldots,\gamma_\ell){(x_1,\ldots,x_\ell)}=\left(\gamma_1{x_1},\ldots,\gamma_\ell{x_\ell}\right)\quad\mbox{and}\quad (\gamma_1,\ldots,\gamma_\ell){\sigma^{-1}}=(\gamma_{1\sigma},\ldots,\gamma_{\ell\sigma})
$$
for all $(\gamma_1,\ldots,\gamma_\ell)\in\Gamma^\ell$,
$x_1,\ldots,x_\ell\in L$, and $\sigma\in H$.
The important properties of wreath products can be found in most
textbooks on permutation group theory, see for instance~\cite{dm}. 

The {\em holomorph} of an abstract 
group $M$ is the semidirect product $M\rtimes\aut M$.
If $M$ is a regular, characteristically simple permutation 
group acting on a set $\Omega$, 
then $\Omega$ can be identified with the underlying set of $M$, and 
$\hol M$ can also be viewed as a subgroup of $\sym\Omega$. 
It is well-known that in this case $\norm{\sym\Omega}M=\hol M$. 

Following~\cite[Section~3]{bad:quasi}, 
we distinguish between 8 classes of finite primitive groups,
namely {\sc HA}, {\sc HS}, {\sc HC}, {\sc SD}, {\sc CD}, {\sc PA}, 
{\sc AS}, {\sc TW}, and 8~classes of finite quasiprimitive groups, namely 
{\sc HA}, {\sc HS}, {\sc HC}, {\sc Sd}, {\sc Cd}, {\sc Pa}, 
{\sc As}, {\sc Tw}. The type of a primitive or quasiprimitive group $G$ can be
recognized from the structure and the 
permutation action of its socle, denoted $\soc G$. Let $G\leq\sym\Omega$ 
be a quasiprimitive 
permutation group, let $M$ be a minimal normal subgroup of $G$, and 
let $\omega\in\Omega$. Note that $M$ is a characteristically simple group, and, if $M$ is non-abelian, then a subdirect 
subgroup of $M$ 
is meant to be subdirect with respect to the unique 
finest direct decomposition of $M$. The
main characteristics of $G$ and $M$ in each primitive and quasiprimitive type
are as follows.

{\sc HA:} $M$ is abelian, $\cent GM=M$ and $G\leq \hol M$. The group $G$ is
always primitive.

{\sc HS:} $M$ is non-abelian, simple, and regular; 
$\soc G=M\times\cent GM$ and $G\leq\hol M$. The group $G$ is always primitive.

{\sc HC:} $M$ is non-abelian, non-simple, and regular; $\soc G=M\times\cent GM$ and $G\leq\hol M$. The group $G$  
is always primitive.

{\sc Sd:} $M$ is non-abelian and non-simple; $M_\omega$ is a 
simple subdirect subgroup of $M$ and $\cent GM=1$. If, 
in addition, $G$ is primitive then the type of $G$ is {\sc SD}.

{\sc Cd:} $M$ is non-abelian and non-simple; $M_\omega$ is a 
non-simple subdirect subgroup of $M$ and  $\cent GM=1$. If, 
in addition, $G$ is primitive then the type of $G$ is {\sc CD}.

{\sc Pa:} $M$ is non-abelian and non-simple; $M_\omega$ is a 
not a subdirect subgroup of $M$ and $M_\omega\neq 1$; $\cent GM=1$. If, 
in addition, $G$ is primitive then the type of $G$ is {\sc PA}.

{\sc As:} $M$ is non-abelian and simple; $\cent GM=1$. If, 
in addition, $G$ is primitive then the type of $G$ is {\sc AS}.

{\sc Tw:} $M$ is non-abelian and non-simple; $M_\omega=1$; $\cent GM=1$. If, 
in addition, $G$ is primitive then the type of $G$ is {\sc TW}.

It is not hard to prove that if 
$G$ is a permutation group
with at least two transitive 
minimal normal subgroups 
then $G$ is primitive of type
{\sc HS} or {\sc HC}.

\section{Cartesian decompositions preserved by quasiprimitive groups}\label{sec2}

A {\em Cartesian decomposition} of a set $\Omega$ is a set
$\{\Gamma_1,\ldots,\Gamma_\ell\}$ of proper partitions of $\Omega$ such that 
$$
|\gamma_1\cap\cdots\cap\gamma_\ell|=1\quad\mbox{for
all}\quad\gamma_1\in\Gamma_1,\ldots,\gamma_\ell\in\Gamma_\ell.
$$
The number $\ell$ is called
the {\em index} of the Cartesian decomposition
$\{\Gamma_1,\ldots,\Gamma_\ell\}$. A Cartesian decomposition is said to be 
{\em homogeneous} if its partitions have the same size.

If $G$ is a permutation group acting on $\Omega$, then a Cartesian
decomposition $\E$ of $\Omega$ is said to be {\em $G$-invariant}, if the partitions in
$\E$ are permuted by $G$. In this case, for $\Gamma\in\E$, 
the permutation group induced by $G_\Gamma$ on $\Gamma$ is denoted
$G^\Gamma$ and is referred to as a {\em component} of $G$.
If $G$ acts transitively on $\E$, then $\E$ is said to be a {\em transitive} $G$-invariant
Cartesian decomposition.

Let $\mathcal E=\{\Gamma_1,\ldots,\Gamma_\ell\}$ be a Cartesian decomposition of
 a set $\Omega$. 
It follows from the last displayed equation that the following map 
is a well-defined bijection between
 $\Omega$ and
$\Gamma_1\times\cdots\times \Gamma_\ell$:
$$
\vartheta:\omega\mapsto(\gamma_1,\ldots,\gamma_\ell)\mbox{ where for }
i=1,\ldots,\ell,\ \gamma_i\in\Gamma_i\mbox{ is chosen so that }
\omega\in\gamma_i.
$$ 
Now suppose that $G$ is a permutation group on $\Omega$ and that $\mathcal E$
is $G$-invariant. Then there is a faithful action of $G$ on $\Gamma_1\times
\cdots\times \Gamma_\ell$ given by 
\[
(\gamma_1,\ldots,\gamma_\ell)g=(\delta_1,\ldots,\delta_\ell)\quad\mbox{for all}\quad
\gamma_1\in\Gamma_1,\ldots,\gamma_\ell\in\Gamma_\ell\quad\mbox{and}\quad g\in G
\]
where $\delta_1\in\Gamma_1,\ldots,\delta_\ell\in\Gamma_\ell$ are defined by
$\{\delta_1,\ldots,\delta_\ell\}=\{\gamma_1g,\ldots,\gamma_\ell g\}$.
(Note that the definition of a Cartesian decomposition ensures that the 
sets $\Gamma_1,\ldots,\Gamma_\ell$ are disjoint.) We observe that 
$(\vartheta,\iota)$, where $\iota:G\to G$ is the identity map on $G$, is a 
permutational isomorphism from $G$ on $\Omega$ to $G$ on $\Gamma_1\times
\cdots\times\Gamma_\ell$, that is 
\[
(\omega\vartheta)g=(\omega g)\vartheta\quad\mbox{for all}\quad \omega\in\Omega\quad\mbox{and}\quad g\in G.
\]
Suppose, in addition, that $\mathcal E$ is homogeneous, and 
set $\Gamma=\Gamma_1$.
Then for $i=1,\ldots,\ell$ there exists a bijection $\alpha_i:\Gamma_i\to\Gamma$
, whence
we have a bijection $\vartheta':\Omega\to\Gamma^\ell$ given by
\begin{equation}\label{sigbij2}
\vartheta':\omega\mapsto(\gamma_1\alpha_1,\ldots,\gamma_\ell\alpha_\ell)
\quad\mbox{for all}\quad \omega\in\Omega,
\end{equation}
where $\omega\vartheta=(\gamma_1,\ldots,\gamma_\ell)$.
Let $\chi$ be the isomorphism $\sym(\Omega)\to\sym(\Gamma^\ell)$
induced by $\vartheta'$.  Then $(\vartheta',\chi)$ restricts to a
permutational isomorphism from $G$ on $\Omega$ to $G\chi$ on $\Gamma^\ell$.
It is clear
that $G\chi$ is contained in the wreath product $\sym(\Gamma)\wr \sy\ell$ 
in its product 
action on
$\Gamma^\ell$. Moreover the image of $\Gamma_i$ under $\vartheta'$ is the
partition of $\Gamma^\ell$ with the parts indexed by $\Gamma$, such
that the $\gamma$-part is the set of all $\ell$-tuples of $\Gamma^\ell$
with $i$-th entry $\gamma$. The set
$\E\vartheta'=\{\Gamma_1\vartheta',\ldots,\Gamma_\ell\vartheta'\}$ is a
$G\chi$-invariant Cartesian decomposition of $\Gamma^\ell$. 

In addition to the assumptions in the previous paragraph, suppose now
$G$ is transitive on $\mathcal E$, and set $\Gamma=\Gamma_1$. 
Then for all $i=1,\ldots,\ell$, there exists an  element $g_i\in G$ such that 
$\Gamma_ig_i=\Gamma$. Define the bijection $\alpha_i:\Gamma_i\to\Gamma$ by 
$\alpha_i:\gamma\mapsto \gamma g_i$ for all $\gamma\in\Gamma_i$. Let $\vartheta'$ be given by 
(\ref{sigbij2}), and, as before, let $\chi:\sym(\Omega)\to\sym(\Gamma^\ell)$ be the 
isomorphism induced by $\vartheta'$. Direct calculation shows that 
$G\chi$ is contained in $G^\Gamma\wr \sy \ell$.
In this paper, we will often identify $G$ with $G\chi$.

Let $M$ be a transitive 
permutation group on a finite set $\Omega$ and let $\E$ be an 
$M$-invariant Cartesian decomposition of $\Omega$. Suppose further that
$M_{(\E)}=M$ and, for $\Gamma\in\E$,  
$M$ can be written as $M=M_1\times M_2$ where
$M_1$ is a normal subgroup of $M$ and $M_2$ is the kernel of the $M$-action on $\Gamma$. It follows that the $M_1$-action on $\Gamma$ must be faithful, and so 
we may naturally identify $M_1$ with $M^\Gamma$. 
In this situation we say that the Cartesian decomposition $\E$ 
is {\em $M$-normal} if
$M=\prod_{\Gamma\in\E}M^\Gamma$. 
If $\E$ is a $G$-invariant Cartesian decomposition of the underlying set of 
some permutation group $G$ then $\E$ is said to be {\em normal} if $\E$ is 
$M$-normal for some transitive normal subgroup $M$ of $G$. In this case 
$M=\prod_{\Gamma\in\E}M^\Gamma$, and, for $\Gamma\in\E$, we denote 
$\prod_{\Gamma_0\neq \Gamma}M^{\Gamma_0}$ by $\overline{M^\Gamma}$. 

In the next lemma we summarise some basic properties of
normal Cartesian decompositions.

\begin{lemma}\label{morbits}
Suppose that $M$ is a transitive  permutation group acting on $\Omega$,
$\E$ is an $M$-normal Cartesian decomposition of $\Omega$,
and let $\Gamma\in\E$.  Then the 
following all hold. 
\begin{itemize}
\item[(a)] The group 
$M^{\Gamma}$ is
transitive on $\Gamma$. 
\item[(b)]  For $\gamma\in\Gamma$ and $\omega\in\gamma$, we have 
$(M^{\Gamma})_{\gamma}=M_\omega\cap M^{\Gamma}$.
\item[(c)] If $\omega\in\Omega$, then
$M_\omega=\prod_{\Gamma'\in\E}\left(M_\omega\cap M^{\Gamma'}\right)$.
\item[(d)] The partition $\Gamma$  is
the set of $\overline{M^{\Gamma}}$-orbits 
and $\overline{M^{\Gamma}}=M_{(\Gamma)}$.
\end{itemize}
\end{lemma}
\begin{proof}
Let $\E=\{\Gamma_1,\ldots,\Gamma_\ell\}$. Without loss of generality we will 
prove parts (a), (b), and (d) in the case $\Gamma=\Gamma_1$.

(a) Suppose 
that $\gamma_1,\ \gamma_2\in\Gamma_1$ and let $\omega_1$ and $\omega_2$ be
arbitrary elements of $\Omega$ such that 
$\omega_1\in\gamma_1$ and $\omega_2\in\gamma_2$.
Then there is some $m\in M$ such that $\omega_1m=\omega_2$. As $\Gamma_1$ 
is a system of imprimitivity for $M$, it follows that 
$\gamma_1m=\gamma_2$. Suppose that $m=m_1\cdots m_\ell$ such that $m_i\in M^{\Gamma_i}$. Then $\gamma_1m=\gamma_1m_1=\gamma_2$. Thus $M^{\Gamma_1}$ is 
transitive on $\Gamma_1$. 

(b) 
Let $\omega\in\Omega$ and 
$\gamma_1\in\Gamma_1$ such that $\omega\in\gamma_1$. 
As the partition $\Gamma_1$ is a system of imprimitivity in $\Omega$ 
for the $M$-action, we obtain  
that $M_\omega\cap M^{\Gamma_1}\leq 
(M^{\Gamma_1})_{\gamma_1}$. Suppose that $m\in (M^{\Gamma_1})_{\gamma_1}$ and that
$\{\omega\}=\gamma_1\cap\gamma_2\cdots\cap\gamma_\ell$ 
for some $\gamma_2\in\Gamma_2,
\ldots,\gamma_\ell\in\Gamma_\ell$. Since $m\in M^{\Gamma_1}$, the element 
$m$ stabilises 
$\gamma_2,\ldots,\gamma_\ell$, and, by assumption, $m$ also stabilises 
$\gamma_1$. Hence $m$ stabilises $\omega$, and so $(M^{\Gamma_1})_{\gamma_1}\leq M_\omega\cap M^{\Gamma_1}$. Therefore  $(M^{\Gamma_1})_{\gamma_1}=M_\omega\cap M^{\Gamma_1}$.

(c)
Assume that 
$\{\omega\}=\gamma_1\cap\cdots\cap\gamma_\ell$ as in part~(b). 
Then 
$(M^{\Gamma_1})_{\gamma_1}
\times\cdots\times (M^{\Gamma_\ell})_{\gamma_\ell}\leq M_\omega$. On the 
other hand, by part~(a), 
\begin{multline*}
|M:(M^{\Gamma_1})_{\gamma_1}
\times\cdots\times (M^{\Gamma_\ell})_{\gamma_\ell}|=|M^{\Gamma_1}:(M^{\Gamma_1})_{\gamma_1}|\times\cdots\times |M^{\Gamma_\ell}:(M^{\Gamma_\ell})_{\gamma_\ell}|\\=
|\Gamma_1|\times\cdots
\times|\Gamma_\ell|=|\Omega|, 
\end{multline*}
and so
$$
M_\omega=(M^{\Gamma_1})_{\gamma_1}
\times\cdots\times (M^{\Gamma_\ell})_{\gamma_\ell}=(M_\omega\cap M^{\Gamma_1})
\times\cdots\times (M_\omega\cap M^{\Gamma_\ell}).
$$

(d) 
Suppose, as above, 
that $\{\omega\}=\gamma_1\cap\cdots\cap\gamma_\ell$ for some $\omega\in\Omega$ and $\gamma_1\in\Gamma_1,\ldots,\gamma_\ell\in\Gamma_\ell$. 
As $\Gamma_1$ is a block-system for the action of $M$ on 
$\Omega$, the block $\gamma_1$ is stabilised by 
$$
K_1=\left(M_\omega\cap M^{\Gamma_1}\right)\times\overline{M^{\Gamma_1}}=(M^{\Gamma_1})_{\gamma_1}\times\overline{M^{\Gamma_1}}. 
$$
On the other hand, by part~(b), $|M:K_1|=|\Gamma_1|$, and so 
$\Gamma_1$ is the system of 
imprimitivity for $M$ 
corresponding
to the overgroup 
$K_1$
of $M_\omega$. Thus $\gamma_1=\omega^{K_1}$. As $\overline{M^{\Gamma_1}}\lhd
M$, the set $\Sigma$ of $\overline{M^{\Gamma_1}}$-orbits is also a system of
imprimitivity for $M$. It suffices to prove that the block
 $\sigma\in\Sigma$ containing $\omega$ is equal to $\gamma_1$. 
If $\omega'\in\gamma_1=\omega^{K_1}$, then there is some $m\in
K_1$, such that $\omega^m=\omega'$. Write $m$ as the product
$m_1\cdots m_\ell$ where $m_i\in M^{\Gamma_i}$ for $i=1,\ldots,\ell$. Then, as
$m_1\in M_\omega$,  we have $\omega^m=\omega^{m_2\cdots m_\ell}$.
Since $m_2\cdots m_\ell\in\overline{M^{\Gamma_1}}$, we obtain
$\omega'\in\sigma$. Therefore $\gamma_1\subseteq\sigma$. On the other
hand, as
$\overline{M^{\Gamma_1}}\leq K_1$, it follows that
$\sigma\subseteq\gamma_1$. Thus $\sigma=\gamma_1$ and the two block
systems $\Sigma$ and $\Gamma_1$ coincide. 

Since $\overline{M^{\Gamma_1}}$ is normal in $M$ and fixes $\gamma_1$, it follows
that $\overline{M^{\Gamma_1}}\leq M_{(\Gamma_1)}$. On the other hand, let 
$m\in M_{(\Gamma_1)}$. Then $m=m_1\cdots m_\ell$, where $m_i\in M^{\Gamma_i}$ for all $i$. As
$m_2\cdots m_\ell\in M_{(\Gamma_1)}$, it follows that $m_1\in
M_{(\Gamma_1)}$ and $m_1$ fixes each $\Gamma_i$ pointwise. Therefore
$m_1$ lies in the kernel of the $M$-action on $\Omega$, and, since $M$ is
faithful, $m_1=1$. This proves that
$M_{(\Gamma_1)}\leq\overline{M^{\Gamma_1}}$, and hence $\overline{M^{\Gamma_1}}=M_{(\Gamma_1)}$. 
\end{proof}

The next result shows that a normal Cartesian decomposition is always normal
with respect to a transitive minimal normal subgroup.

\begin{theorem}\label{Unilem}
Let $G$ be a permutation group on $\Omega$, 
$M$ a transitive, non-abelian, minimal normal subgroup of $G$, and let 
$\E$ be
a normal $G$-invariant Cartesian decomposition of $\Omega$.
Then $\E$ is $M$-normal  and $G$ is transitive on $\E$.
\end{theorem}
\begin{proof}
Let $N$ be a normal subgroup of $G$, such that $\E$ is $N$-normal, and
let $N_i=N^{\Gamma_i}$, where $\E=\{\Gamma_1,\ldots,\Gamma_\ell\}$. By the definition of 
$N$-normal, $N=N_1\times\dots\times N_\ell$.
Suppose first that $M\leq N$. We claim that 
\begin{equation}\label{dirdec}
M=(N_1\cap M)\times\cdots\times (N_\ell\cap M).
\end{equation}
Note that $N=N_i\times \overline{N_i}$ and $M\trianglelefteq N$. Since
$M$ is a non-abelian minimal normal subgroup of $G$, $M$ is a direct
product of isomorphic non-abelian simple groups. Let $T$ be a simple direct factor
of $M$. For $i\in \{1,\ldots,\ell\}$ let $\sigma_{N_i}$ be the projection map
$N\to N_i$. As $T$ is a non-trivial subgroup of $N$, there exists $i$ such that $\sigma_{N_i}(T)$ is
non-trivial, whence $\sigma_{N_i}(T)\cong T$ as $T$ is a non-abelian simple group. We
 claim that
$T\leq N_i$.
Choose $x\in \sigma_{N_i}(T)$ with $x\not=1$. As $M$ is a normal subgroup of 
$G$, the subgroup $T^x$ is also a minimal normal subgroup of $M$, and we 
see that either
$T=T^x$ or $[T,T^x]$ is trivial. If the latter, then 
\[
1=\sigma_{N_i}([T,T^x])=[\sigma_{N_i}(T),(\sigma_{N_i}(T))^x]=[\sigma_{N_i}(T),\sigma
_{N_i}(T)];
\]
but the last term is non-trivial as $\sigma_{N_i}(T)$ is  a  
non-abelian simple group. Hence the former holds, that is $T=T^x$.
As $\sigma_{N_i}(T)\cong T$, if $t,\,t'\in T$ then
$t^x=t'$ if and only if $\sigma_{N_i}(t)^x=\sigma_{N_i}(t')$.
As $\sigma_{N_i}(T)$ is non-abelian simple, conjugation by $x$ induces a 
non-trivial automorphism of $\sigma_{N_i}(T)$. Thus it must induce a non-trivial
automorphism of $T$ as well.
However if $\sigma_{N_j}(T)$ is non-trivial for some $j\neq i$, 
then $x$ centralises $\sigma_{N_j}(T)\cong T$ and so, as shown by a similar argument, conjugation by $x$ induces a trivial automorphism of $T$. Thus $\sigma_{N_j}(T)$ is
 trivial for all $j\ne i$, and
$T\leq N_i$ as claimed. Thus each simple direct factor of $M$ is
contained in some $N_i$ and it follows that~\eqref{dirdec} holds.
In this case, therefore, $\E$ is
$M$-normal, so $M=\prod_{\Gamma\in\E}M^\Gamma$. Since $G$ is transitive on the simple
 direct factors of
$M$, it follows that $G$ is transitive on $\E$. 

Thus we may assume that
$M\not\leq N$, so $M\cap N=1$, by minimality of $M$. As both $M$ and
$N$ are transitive, we have that they are both regular and, for a fixed
$\omega\in\Omega$, the map
$\vartheta:N\rightarrow M$ given by 
$$
x\vartheta=y^{-1}\quad\mbox{if and only if}\quad \omega x=\omega y
$$
is an isomorphism between $N$ and $M$. For each $i$ set
$M_i=N_i\vartheta$ and $\overline{M_i}=\overline{N_i}\vartheta$. 
Since $\vartheta$ is an
isomorphism, $M=M_1\times\cdots\times M_\ell$. Also by the
definition of $\vartheta$, the $\overline{M_i}$-orbits are the same as
the $\overline{N_i}$-orbits, and by Lemma~\ref{morbits}(d) the
$\overline{M_i}$-orbits form the partition $\Gamma_i$. Thus
$M^{\Gamma_i}=(M_i)^{\Gamma_i}$. If $K_i$ is the kernel of the action of
$M_i$ on $\Gamma_i$, then $K_i\vartheta^{-1}$ (by the definition of
$\vartheta$) also acts trivially on $\Gamma_i$, and hence lies in
$N_{(\Gamma_i)}=\overline{N_i}$. However $K_i\vartheta^{-1}\leq
N_i$, so $K_i\vartheta^{-1}\leq N_i\cap\overline{N_i}=1$, whence
also $K_i=1$. Hence $M^{\Gamma_i}=M_i$, and $M=\prod_{i=1}^\ell
M^{\Gamma_i}$, so $\E$ is $M$-normal.
As $M$ is a minimal normal subgroup of $G$,
the conjugation action by $G$ on the $M^{\Gamma_i}$ is transitive, and
hence $G$ is transitive on $\E$.
\end{proof}

Motivated 
partly by Theorem~\ref{Unilem},
we now introduce the class of innately transitive permutation groups.
A finite permutation group is said to be {\em innately transitive} if it has 
a transitive minimal normal subgroup (see~\cite{bp} for a comprehensive study
of these groups). In particular, primitive and quasiprimitive
groups are innately transitive. In~\cite{transcs} we introduced six disjoint classes
of transitive Cartesian decompositions that may be preserved by an
innately transitive group. 
Theorem~\ref{Unilem} has the following consequence: if $\E$
is a transitive $G$-invariant normal Cartesian decomposition,
for an innately transitive group $G$, then $\E$ belongs to one
of only two of the six classes in~\cite{transcs}, namely 
$\cd{S}G$ or $\cd 1G$ (see~\cite{transcs} for 
the notation). Moreover, the Cartesian decompositions in these
two families are normal.
This simple 
observation will be refined somewhat  in Theorem~\ref{main}(d).


Suppose that $\E$ is a $G$-invariant 
Cartesian decomposition for some permutation group $G$.
We say that ${\mathcal E}$ is a  {\em blow-up
  decomposition} for $G$ if ${\mathcal E}$ is 
transitive and it is $M$-normal for some transitive 
normal subgroup $M$ of $G$ 
such that, for all  $\Gamma\in\E$,
we have  $M^{\Gamma}=\soc(G^{\Gamma})$. 



The concept of a `blow-up' is due to Kov\'acs \cite{blowups}. In the above we have simply
translated his definition to the current context. 
The terminology `blow-up' is intended to stress the
intuitive idea that a permutation group $G$ on a set $\Omega$ with a
blow-up decomposition $\mathcal E$ is simply a `blown-up'  version of the 
smaller permutation group $G^\Gamma$ for $\Gamma\in\mathcal E$; thus we talk of $G$
as being a `blow-up' of its components.
Likewise, if $\mathcal E$ is a Cartesian decomposition, then there is a sense
in which $G$ can be thought of as being constructed from its components (which are necessarily 
groups smaller than $G$), although in general the relationship between $G$ and its components is not 
as strong as in the blow-up case.

\section{Components and blow-up decompositions in quasiprimitive groups}\label{sec3}

The aim of this section is to study the relationship between
a quasiprimitive permutation group and its components 
under a blow-up decomposition.

\begin{lemma}\label{5.3}
Let  $G$ be a permutation group on a set $\Omega$ and suppose  that 
$\mathcal E$ is a blow-up decomposition for $G$
such that, for $\Gamma\in \mathcal
E$, $G^\Gamma$ is quasiprimitive on $\Gamma$.
Suppose further that  either $\soc(G^\Gamma)$ is non-abelian or $G$ is quasiprimitive. 
Then the rule $K\mapsto 
K^\Gamma$ defines a bijection from the set of minimal 
normal subgroups of $G$ to the set of minimal normal subgroups of $G^\Gamma$.
\end{lemma}

\begin{proof}
To prove this we adapt the argument preceding~\cite[(2.1)]{blowups}. Set $H=G^\Gamma$, $M=\soc H$,  and $\ell =|\mathcal E|$. As explained in Section~\ref{sec2}, 
we may assume without loss of generality that $G$ is a subgroup of the wreath product $W=H\wr \sy\ell$ in 
its product action on $\Gamma^\ell$, and (since $\E$ is a blow-up
decomposition)
that  $G$ contains $M^\ell=(\soc H)^\ell$.
Since $G^\Gamma$ is quasiprimitive, $\cent HM\leq M$. It follows easily that 
$\cent {W}{M^\ell}\leq M^\ell$.
Thus each minimal normal subgroup of $G$ must lie in 
$M^\ell$, and hence must lie in $G_{(\E)}$. 

If $M$ is non-abelian, then $M$, 
and also $M^\ell$, 
is a direct 
product of non-abelian simple groups. 
Given  that $M^\ell\leq G$, we deduce that 
each minimal normal subgroup of $G$ is  a direct product
of the $G$-conjugates of some simple direct factor of $M^\ell$. 
Since $G$ is transitive on $\mathcal E$, the projection of $G$ onto the 
top group $\sy\ell$ of the wreath product $W$ is a transitive subgroup of
$\sy\ell$. Thus each minimal normal subgroup $K$ of $G$ is of the form $K_0^\ell$, where $K_0$ is a characteristically simple normal subgroup of $H$. Moreover, the minimality of $K$ implies that $H=G^\Gamma$ is transitive on the simple direct factors of $K_0$, and so $K_0=K^\Gamma$ is a minimal normal subgroup of $H$. Conversely for each minimal normal subgroup $K_0$ of $H$, $K_0^\ell$ is a minimal normal subgroup of $G$.

If $M$ is abelian then, by assumption, in this case, $G$ is quasiprimitive. 
Moreover $G$ has a minimal normal subgroup contained in $M^\ell$ 
that is abelian, and hence $G$ is quasiprimitive of type HA. 
This implies that $G$ is primitive.
We now apply~\cite[(2.1)]{blowups} directly
to deduce that $\soc G=M^\ell$. As both $G$ and $H$ are now quasiprimitive with an abelian socle, 
they both have a unique minimal normal subgroup, namely $\soc G$ and $M$
respectively. This gives the required result.
\end{proof}

\begin{corollary}\label{5.35}
Suppose that $G$ is a permutation group on $\Omega$ and  $\E$ is a blow-up
decomposition for $G$. Let $\Gamma_0\in\E$. 
\begin{enumerate}
\item[(i)] If $G$ is quasiprimitive then the component $G^{\Gamma_0}$ is 
quasiprimitive and $\soc G=\prod_{\Gamma\in\E}\soc(G^{\Gamma})$.
\item[(ii)] If the component  $G^{\Gamma_0}$ is  quasiprimitive and 
$\soc (G^{\Gamma_0})$ is non-abelian, then  $G$ is quasiprimitive.
\end{enumerate}
\end{corollary}
\begin{proof}
(i) 
Set $H=G^{\Gamma_0}$ and $\ell =|\mathcal E|$.
We may assume that $G$ is a quasiprimitive subgroup of the wreath product $H\wr \sy\ell$ in 
its product action on $(\Gamma_0)^\ell$. Since, $\E$ is a blow-up decomposition, 
$G$ contains $(\soc H)^\ell$.
If $H$ is not quasiprimitive, then there
exists a minimal normal subgroup $K$ of $H$ with $K$ intransitive on 
$\Gamma_0$. 
Then $K^\ell\leq (\soc H)^\ell$ and $K^\ell$ is a normal subgroup of 
$H\wr \sy\ell$ contained in $G$ and is 
intransitive on $(\Gamma_0)^\ell$. This contradicts the quasiprimitivity of $G$.
Hence $H$ is quasiprimitive.
Then by Lemma~\ref{5.3}, for a minimal normal subgroup $K$  of $G$, 
$K^{\Gamma_0}$ is a minimal normal subgroup of $G^{\Gamma_0}$ which 
implies the assertion.

(ii) By Lemma~\ref{5.3}, a minimal normal subgroup of $G$
is of the form $\prod_{\Gamma\in\E} K^\Gamma$ where, for each
$\Gamma$, $K^\Gamma$ is a
minimal normal subgroup of $G^\Gamma$. Since, for $\Gamma\in\E$,  the group
$G^{\Gamma}$ is
quasiprimitive, $K^{\Gamma}$ is transitive on $\Gamma$, and hence
$\prod_{\Gamma\in\E}K^\Gamma$ is transitive on $\Omega$. 
\end{proof}

We interpret the previous results as saying that, given a group 
$G$ with a blow-up decomposition,
the structure of the socle of $G$ is strongly related to the 
structure of the socles of the  components of $G$. 
Also there is a strong link 
between possible quasiprimitivity of $G$ and of its components. 
The following 
example shows the necessity of the restriction that 
$\soc(G^{\Gamma_0})$ cannot be abelian in Corollary~\ref{5.35}(ii).

\begin{example}
Suppose that $H=\left<(1,2,3)\right>$ acting on $\Gamma=\{1,2,3\}$ and let
$G=H\wr D_8$ where $D_8$ is the dihedral group acting on the set $\{1,2,3,4\}$
preserving the block system $\{\{1,2\},\{3,4\}\}$. We consider $G$ as a 
permutation group acting in product action on $\Gamma^4$. Then $G^\Gamma=H$,
and so $G^\Gamma$ is quasiprimitive with a unique minimal normal subgroup 
(namely itself). On the other hand 
$G$ has three minimal normal subgroups 
\begin{eqnarray*}
M_1&=&\{(x,x,x,x)\ |\ x\in H\};\\
M_2&=&\{(x,x,x^2,x^2)\ |\ x\in H\};\\
M_3&=&\{(x,x^2,y,y^2)\ |\ x,\ y\in H\}.
\end{eqnarray*}
and $(M_1)^\Gamma=(M_2)^\Gamma=(M_3)^\Gamma=H$.
So $G$ is not quasiprimitive, and hence the condition in Corollary~\ref{5.35}(ii)
that $\soc(G^{\Gamma_0})$ is non-abelian is necessary. This example also shows that the correspondence
$M\mapsto M^\Gamma$ in Lemma~\ref{5.3} is not always one-to-one if $\soc(G^\Gamma)$ is abelian but $G$ is not quasiprimitive.
\end{example}

The statement of the primitive analogue of Corollary~\ref{5.35}, 
can be obtained by replacing all occurrences of ``quasiprimitive'' by ``primitive'', and the single occurrence of 
``non-abelian'' by ``non-regular'' in 
the statements of these theorems. The validity of these  analogues follows from
\cite[Theorem~1 and (2.1)]{blowups}. Given that for primitive groups
the restriction to non-regular socle
is a stronger condition than the restriction to non-abelian socle, 
we see that the concept of a blow-up
in fact behaves better with respect to quasiprimitivity than it does to primitivity.

\section{Normal decompositions and blow-up decompositions}\label{sec4}

The observant reader may ask whether, for a permutation group $G$,
it is possible that a transitive, 
$G$-invariant, normal decomposition  with
quasiprimitive components is not a blow-up. 
This question is answered by the results of this section.

Suppose that $G\leq\sym\Omega$ is a finite permutation group
with a non-abelian, non-simple, regular,
minimal normal subgroup $M_1$. 
The centraliser $\cent{\sym\Omega}{M_1}$
is isomorphic to $M_1$, and so it is isomorphic to $T^k$ where $T$ is a 
non-abelian finite simple group, and $k\geq 2$. 
If $\cent G{M_1}$
is a proper subdirect subgroup of $\cent{\sym\Omega}{M_1}$, then,
using the terminology of~\cite{bp}, 
$G$ is said to be an  {\em innately transitive group of diagonal 
quotient type}. 
In this case $G$ has
two minimal normal subgroups $M_1$ and $N_1$
where $M_1\cong T^k$ and $N_1\cong T^{k/m}$ for some 
divisor $m$ of $k$ such that $1<m\leq k$.
Further, $N_1$ is semiregular and intransitive, and, in particular, $G$ is
not quasiprimitive.

\begin{theorem}\label{r2}
Let $G$ be a permutation group on a 
set $\Omega$  and suppose that $\mathcal E$ is a transitive,
normal, $G$-invariant Cartesian decomposition of $\Omega$, and that,
for $\Gamma\in \mathcal E$, the component  
$G^\Gamma$ is quasiprimitive. Then  exactly one of the following possibilities holds.
\begin{enumerate}
\item $\mathcal E$ is a blow-up decomposition.
\item $G$ is quasiprimitive of type {\sc Tw} and $G^\Gamma$ 
is primitive of type {\sc HS} 
or {\sc HC}.
\item $G$ is innately transitive with diagonal quotient type and $G^\Gamma$ is primitive of type {\sc HC} or 
{\sc HS}.
\end{enumerate}
\end{theorem}
\begin{proof}
Set $H=G^\Gamma$ and $\ell=|\mathcal E|$. We may assume that $G$ is a subgroup of 
$H\wr \sy\ell$ in its product action on $\Gamma^\ell$.
As $\mathcal E$ is a
normal Cartesian decomposition, there exists a transitive normal subgroup 
$M$ of $G$ such that $M\leq G_{(\mathcal E)}\leq H^\ell$
and $M=(M^\Gamma)^\ell$. 
If $\soc H\leq M^\Gamma$, then by the transitivity of $G$ on $\mathcal E$, we 
have $(\soc H)^\ell\leq G$. Hence $\mathcal E$ is a blow-up decomposition and
part~(i) holds. 
Assume now that 
$\soc H\not\leq M^\Gamma$. As $M^\Gamma$ is normal in $H$, it follows
that $H$ has at least two minimal normal subgroups.
As discussed in Section~\ref{newsec}, the group $H$ is primitive of type HS or HC, and $H$ has exactly two minimal normal subgroups $M_1$, $N_1$,
both regular, non-abelian and transitive on $\Gamma$. Without loss of generality
 we may 
assume that  $M_1\leq M^\Gamma$. Then $(M_1)^\ell\leq G$ and therefore, by definition, $\mathcal
 E$ is $(M_1)^\ell$-normal, so we may assume that $M=(M_1)^\ell$. 
Thus $M$ is a regular non-abelian minimal normal subgroup of $G$. Let 
$C=\cent{\sym\Omega}{M}$ and note that $\soc G\leq M\times C$. 
By~\cite[Theorem~4.3B]{dm}, $C\cong M$. Since $N_1=\cent{\sym\Gamma}
{M_1}\leq H$,
we have $C=(N_1)^\ell\leq H\wr \sy\ell$; therefore
$\soc(H \wr \sy\ell) = M\times C$. We may write $C=\prod_{s\in S} T_s$, where $S
$ is a set of size $k\ell$, and each $T_s$ is isomorphic to a non-abelian simple
 group $T$. 
As $H=G^\Gamma$ is transitive on the set of minimal normal subgroups of 
$N_1=C^\Gamma$ and $G$ is transitive on $\E$, we have that 
$G$ induces a transitive permutation group on $S$ of degree $k\ell$.
If $C\leq G$ then $\soc G=M\times C$, and so $G$ is quasiprimitive 
with two minimal normal subgroups.
In this case $\soc(G^\Gamma)=M_1\times N_1$ and 
$(\soc (G^\Gamma))^\ell=M\times C\leq G$. Therefore $\E$ is a blow-up
decomposition and part~(i) holds.
If $C\cap G=1$ then $G$ 
has a unique minimal normal subgroup, which is regular. Therefore $G$ is
quasiprimitive of type {\sc Tw} and part~(ii) holds.

Thus we may assume that $1< C\cap G< C$. 
In this case $1\neq (C\cap G)^\Gamma\lhd G^\Gamma$ and 
$(C\cap G)^\Gamma\leq C^\Gamma$. Since $G^\Gamma$ is quasiprimitive, 
$(C\cap G)^\Gamma$ must be transitive. Recall that $C^\Gamma$ is
regular. Thus
$(C\cap G)^\Gamma=C^\Gamma\cong T^k$, and also $N=C\cap G=\cent GM\neq 1$. 
Therefore,
$N$ is a proper subdirect subgroup of $C$ where $C$ is viewed as a
direct product of its minimal normal subgroups, and hence $N$ is a
direct product of full diagonal subgroups. 
Thus $G$ is innately transitive with diagonal quotient type.

These possibilities are mutually exclusive. For if $\E$ is a 
blow-up and $\soc G$ is non-abelian then the other two possibilities cannot 
occur, by Corollary~\ref{5.35}. On the other hand, if $\soc G$ is 
abelian then only possibility~(i) can occur.
\end{proof}

We construct an example to show that the situation described by
Theorem~\ref{r2}(iii) is possible.

\begin{example}
Let $T$ be a non-abelian finite simple group, 
let $H$ be any subgroup
of the holomorph 
$$
{\rm Hol}(T^k)=T^k\rtimes\aut(T^k)=T^k\rtimes({\rm Aut}(T) \wr \sy k)
$$
such that $H$ has two minimal normal subgroups $M_1$ and $N_1$
where $M_1\cong N_1\cong T^k$. Then $H$, considered
as a permutation group on $\Gamma=T^k$,  is a primitive group of type
HS if $k=1$ or type HC if $k\geq 2$. 
Let $\ell\geq 2$, and let $G$ be a subgroup of $H \wr \sy \ell$ in
its product action on $\Omega=\Gamma^\ell$, such that $G$ contains $M_1^\ell\cong
(T^k)^\ell$,
a regular normal subgroup, and $G = M_1^\ell(H\delta \times
\sy \ell$) where $\delta$ is  the diagonal embedding 
$\delta:H\rightarrow H^\ell$ defined by $h\delta=(h,\ldots,h)$. 
Then $\soc G = M_1^\ell \times N_1\delta$, and, in addition, 
$N_1\delta$  is a semiregular and intransitive
minimal normal subgroup of $G$. Thus $G$ is innately transitive
with diagonal quotient type in its action on $\Gamma^\ell$. Moreover, the component of $G$ induced on $\Gamma$ is $H$,
which is primitive. 
\end{example}

Theorem~\ref{r2} demonstrates that the quasiprimitivity of the components
of a transitive normal Cartesian decomposition for $G$ often implies
that $G$ also is quasiprimitive. More precisely, the following result
is valid.

\begin{corollary}
Suppose that $G\leq\sym\Omega$ is a permutation group and $\E$ is a 
transitive, 
normal $G$-invariant Cartesian decomposition of $\Omega$ such that, 
for $\Gamma\in\E$, the component $G^\Gamma$ is quasiprimitive but not
of type {\sc HA}, {\sc HS} or {\sc HC}.  Then $\E$ is a blow-up
decomposition and $G$ is a quasiprimitive group.
\end{corollary}
 \begin{proof}
As the type of a component is not HS or HC, 
it follows from Theorem~\ref{r2} that $\E$ is a blow-up decomposition. Then Corollary~\ref{5.35}(ii) 
implies that 
$G$ is quasiprimitive.
\end{proof}

Next in this section we present a sufficient and necessary condition
to decide if a transitive Cartesian decomposition is a blow-up 
decomposition.

\begin{theorem}\label{r2.5}
Let  $G$ be a quasiprimitive permutation group  acting on a set $\Omega$ and let  
$\mathcal E$ be a transitive $G$-invariant Cartesian decomposition of 
$\Omega$.
Then $\mathcal E$ is a blow-up decomposition if and only if $\E$ is $(\soc G)$-normal and, for $\Gamma\in\E$, $\cent{G^\Gamma}{(\soc G)^\Gamma}\leq (\soc G)^\Gamma$.
\end{theorem}

\begin{proof}
Set $H=G^\Gamma$, $\ell=|\mathcal E|$,  and $M=\soc G$. We may assume that $G$ is a subgroup
of $H\wr \sy \ell$ in its product action on $\Gamma^\ell$.
If $\E$ is  a blow-up decomposition, then, 
by Corollary~\ref{5.35}(i), $M=(\soc H)^\ell$. In particular,
$\E$ is $M$-normal and $M^\Gamma=\soc H$. Furthermore, 
Corollary~\ref{5.35}(i) implies that 
the component  $H$ is quasiprimitive whence we have
$\cent H{M^\Gamma}\leq M^\Gamma$
as required. 
 
Conversely, assume that $\mathcal E$ is $M$-normal, and that 
the centraliser
$\cent H{M^\Gamma}$ is contained in $M^\Gamma$. By the definition of a blow-up it is enough to show
that $M^\Gamma=\soc H$. 
If $M$ is abelian then $G$ is of type 
HA and $M$ is the unique minimal normal subgroup of $G$. We also know 
that $M$ is regular and elementary 
abelian, hence so is $M^\Gamma$. As $M^\Gamma\unlhd H$ we must have 
that $H\leq {\sf Hol}(M^\Gamma)$, the normaliser of $M^\Gamma$ 
in $\sym\Gamma$. The subgroup $M^\Gamma$ must be a minimal normal 
subgroup of $H$, for if $N<M^\Gamma$ and $N\unlhd H$, then 
$N^\ell<M$ is a normal subgroup of $G$ properly contained in $M$, which is 
impossible by the minimality of $M$. 
Therefore $H$ is primitive of type HA on $\Gamma$ and $M^\Gamma=\soc H$.

Thus we may assume that 
$M$ is non-abelian. Then $M^\Gamma$, which is a homomorphic image of $M$, is the 
direct product of non-abelian simple groups and so $M^\Gamma$ is  necessarily a
direct product of  minimal normal subgroups of $H$. However, as $\cent
H{M^\Gamma}\leq M^\Gamma$ 
it follows that $M^\Gamma$ must contain all minimal normal subgroups of $H$, whence $M^\Gamma=\soc H$ as required.
\end{proof}

We end this section with an example to show that a quasiprimitive group
may have non-quasiprimitive components with respect to a normal 
Cartesian decomposition.

\begin{example}\label{ex-1}
Let $T$ be a non-abelian finite simple group and let $P$ be a finite group
with a core-free 
subgroup $Q$ and a homomorphism $\varphi:Q\rightarrow\aut T$ such that
$\varphi$ induces a non-trivial, proper subgroup $R$ of $\inn T$. Then the 
twisted wreath product $W=T\wr_\varphi P$ is a quasiprimitive permutation
group acting on $\Omega=T^{|P:Q|}$; see~\cite{tw}. Further, the
natural Cartesian decomposition of $\Omega$ is clearly $(\soc W)$-normal.
However, a component of $W$ has a 
regular minimal normal subgroup isomorphic to
$T$, and an intransitive normal subgroup
isomorphic to $R$. Thus a component of $W$ is not quasiprimitive.
\end{example}

\section{Inclusions of quasiprimitive groups}\label{sec5}

In~\cite{bad:quasi} the first two authors studied inclusions of quasiprimitive
groups into primitive ones. The description of such inclusions in the case 
when the primitive group has type {\sc PA} relied on the following theorem,
which was stated without proof in~\cite[Theorem~4.7]{bad:quasi}. Here we give the first published proof. Recall that a Cartesian decomposition $\E$ is homogeneous if $|\Gamma|$ is the same for all $\Gamma\in\E$.

\begin{theorem}\label{main}
If $G\leq\sym\Omega$ 
is a transitive permutation group, $\E$ is
a homogeneous $G$-invariant Cartesian decomposition of $\Omega$, 
and $\Gamma_0\in\E$, 
then the following all hold.
\begin{enumerate}
\item[(a)] If $\E$ is a blow-up decomposition and $G$ is quasiprimitive on 
$\Omega$ then the component $G^{\Gamma_0}$ is quasiprimitive and 
$\soc G=
\prod_{\Gamma\in\E}\soc(G^\Gamma)$.
\item[(b)] If $\E$ is a blow-up decomposition and the component $G^{\Gamma_0}$ is 
quasiprimitive on $\Gamma_0$ not of type {\rm HA} then 
$G$ is quasiprimitive on $\Omega$ and $\soc G=
\prod_{\Gamma\in\E}\soc(G^\Gamma)$.
\item[(c)] The group $G$ is not quasiprimitive of type {\sc Sd}.
\item[(d)] If $G$ is quasiprimitive of type {\sc Cd} then $\E$ is a blow-up 
decomposition and the component $G^{\Gamma_0}$ 
is quasiprimitive of type {\sc Sd} or
{\sc Cd}.
\end{enumerate}
\end{theorem}
\begin{proof}
Parts~(a) and~(b) follow from Corollary~\ref{5.35}. By~\cite[Corollary~1.3]{transcs}, 
quasiprimitive groups of type {\sc Sd} do not preserve Cartesian 
decompositions, and so part~(c) also holds.

Let us now prove part~(d). Suppose that $G$ is a quasiprimitive group of type
{\sc Cd}.
Let $M$ denote the socle of $G$. Then
$M$ is a minimal normal subgroup of $G$. Hence $M$ is a non-abelian characteristically simple group and so 
$M=T^k$ where $T$ is a  non-abelian finite simple 
group. 
Moreover, a point stabiliser $M_\omega$ is a subdirect subgroup of $M$.
It follows from~\cite[Theorem~1.2(c)]{transcs} that $\E$ is $M$-normal
and that $G$ is transitive on $\E$. 
Thus $M$ may be written as 
$M=\prod_{\Gamma\in\E}M^{\Gamma}$. Further, for $\Gamma\in\E$, 
the subgroups $M^{\Gamma}$  are permuted by
$G$. By Lemma~\ref{morbits}, 
the partition $\Gamma$
is the set of $\overline{M^{\Gamma}}$-orbits of 
$\Omega$, and so 
$G_{\Gamma}=\norm{G}{\overline{M^{\Gamma}}}=\norm G{M^{\Gamma}}$.
The subgroup
$M^{\Gamma}$ 
is a transitive minimal normal subgroup of $G^{\Gamma}$. 
As a point stabiliser
in $M^{\Gamma}$ is a subdirect subgroup, 
we obtain from~\cite[Proposition~5.5]{bp} that
$G^{\Gamma}$ is quasiprimitive of type {\sc Sd} or {\sc Cd}. Therefore
Theorem~\ref{r2} implies that $\E$ is a blow-up decomposition.
\end{proof}

\end{document}